\documentclass{amsart}
\usepackage{amssymb,amscd}

\newcounter{intro}

\newtheorem{theo}[intro]{Theorem}

\newtheorem{thm}{Theorem}[section]
\newtheorem{lem}[thm]{Lemma}
\newtheorem{prop}[thm]{Proposition}
\newtheorem{cor}[thm]{Corollary}

\theoremstyle{remark}
\newtheorem{rem}[thm]{Remark}

\numberwithin{equation}{section}

\numberwithin{equation}{section}
\newcounter{counteroman}
\newenvironment{enumeroman}{\begin{list}{\roman{counteroman})}{\usecounter{counteroman}}}{\end{list}}

\newcommand{\cref}[1]{Corollary~\ref{#1}}

\newcommand{\R}{\mathbb{R}}

\newcommand{\cS}{\mathcal{S}}

\let\De=\Delta
\let\eps=\varepsilon

\DeclareMathOperator{\supp}{supp}

\DeclareMathOperator{\vol}{vol}

\def\cro#1#2{\mathrel{\langle {#1},{#2}\rangle}}

\def\Ci{C_0^\infty(M)}
\def\eucl{\rm eucl}
\begin{document}

\title[Riesz transforms on connected sums]
{Riesz transforms on connected sums}
\author{Gilles Carron}
\address{Laboratoire de Math\'ematiques Jean Leray (UMR 6629), 
Universit\'e de Nantes,
2, rue de la Houssini\`ere, B.P.~92208, 44322 Nantes Cedex~3, France}
\email{Gilles.Carron@math.univ-nantes.fr}
\maketitle
\section{Introduction}
Let $(M,g)$ be a complete Riemannian manifold with infinite volume, we denote by $\De=\De^g$ its Laplace
operator, it has an unique self-adjoint extension on $L^2(M,d\vol_g)$ which is also denoted
by $\De$. The Green formula and the spectral theorem show that  for any $\varphi\in \Ci$ :
$$\|d\varphi\|^2_{L^2}=\cro{\De\varphi}{\varphi}=\|\De^{1/2}\varphi\|^2_{L^2};$$
hence the Riesz transform $T:=d\De^{-1/2}$ extends to a bounded operator 
$$T\,:L^2(M)\rightarrow L^2(M; TM).$$

On the Euclidean space, it is well known that the Riesz transform has also a bounded extension 
$L^p(M) \to L^p(M; TM)$ for any $p\in ]1,\infty[$. However, this is not a general feature of 
the Riesz transform on complete Riemannian manifolds, as the matter of fact,
 on the connected sum of two copies of the Euclidean space $\R^n$, the Riesz transform is not bounded on
 $L^p$ for any $p\in [n,\infty[\cap ]2,\infty[$ (\cite{CD, CCH}). It is of interest
  to figure out the range of $p$ for  which $T$ extends to a bounded map $L^p(M) \to L^p(M; T^*M)$. The main
  result of \cite{CCH} answered to this question for manifolds with Euclidean ends :
  
\begin{theo}\label{thA}
 Let $M$ be a complete Riemannian manifold of dimension
  $n \geq 3$ which is the union of a compact part and a finite number of Euclidean ends. 
  Then the Riesz transform is bounded from $L^p(M)$ to $L^p(M; T^*M)$ for $1 < p < n$,
   and is unbounded on $L^p$ for all other values of $p$
if the number of ends is at least two. 
\end{theo}

The proof of this result used an asymptotic expansion of the Schwarz kernel of the resolvent 
$(\De+k^2)^{-1}$ near $k\to 0$. In \cite{CCH} using $L^p$ cohomology, we also find a criterion
 which insures that the Riesz transform
is unbounded on $L^p$ :
\begin{theo}\label{conex} Assume that $(M,g)$ is a complete Riemannian manifold with Ricci curvature bounded
 from below such that
for some $\nu>2$ and $C>0$ $(M,g)$ satisfies the Sobolev inequality
$$\forall \varphi\in  \Ci, C\|\varphi\|_{L^{\frac{2\nu}{\nu-2}}}\le \|d\varphi\|_{L^2}$$
and 
\begin{equation}
\label{volumgro}\forall x\in M
, \forall r>1,\ \vol B(x,r)\le C r^\nu.
\end{equation}
If $M$ has at least two ends, then the Riesz transform is not bounded on $L^p$ for any $p\ge \nu$
\end{theo}

Let $(N,g_0)$ be a simply connected nilpotent Lie group of dimension $n>2$
(endowed with a left invariant metric).  According to \cite{Alex} 
we know that the Riesz transform on $(N,g_0)$ is bounded on $L^p$
for every $p\in ]1,\infty[$.  Let
$\nu$ be the homogeneous dimension of $N$;  for instance we can set
$$\nu=\lim_{R\to \infty} \frac{\log \vol B(o,R)}{\log R},$$
$o\in N$ being a fixed point. Let $(M,g)$ be a manifold isometric at infinity to $k>1$ copies of
$(N,g_0)$. That is to say there are compact sets $K\subset M$ and $K_0\subset N$ such that
$(M\setminus K, g)$ is isometric to $k$ copies of $(N\setminus K_0, g_0)$.
According to \cite{CD} we know that on $(M,g)$ the Riesz transform is bounded on $L^p$ for $p\in ]1,2]$. 
 And the theorem \ref{conex} says that  the Riesz transform is not bounded 
on $L^p$ when $p\ge \nu$. In \cite{CCH}, we make the following conjecture :
{\it Show that the Riesz transform on $(M,g)$ is bounded on $L^p$ for $p\in]1,\nu[$.\\}
The main result of this paper gives a positive answer to this conjecture ; in fact we obtain a more general
result concerning the boundedness of Riesz transform for connected sums, under some mild geometrical conditions :
\begin{theo}\label{main}
Let $(M_0,g_0)$ be a complete Riemannian manifold, we assume that the Ricci curvature of $(M_0,g_0)$
is bounded from below, that the injectivity radius of $(M_0,g_0)$ is positive 
and that for some $\nu>3$ and $C>0$,  $(M_0,g_0)$ satisfies the Sobolev inequality
$$\forall \varphi\in  C_0^\infty(M_0),\ \, C\|\varphi\|_{L^{\frac{2\nu}{\nu-2}}}\le \|d\varphi\|_{L^2}$$
If on $(M_0,g_0)$ the Riesz transform is bounded on $L^p$ for some $p\in ]\nu/(\nu-1), \nu[$, then the
Riesz transform is also bounded on $L^p$ for any manifold $M$ isometric at infinity to several copies of
$(M_0,g_0)$.
\end{theo}

Moreover under a uniform upper growth control of the volume of geodesic balls (such as (\ref{volumgro})), 
the result of \cite{CD} implies that the Riesz transform is bounded on $M$ for any $p\in ]1,2]$ ; hence the
restriction of $p>\nu/(\nu-1)$ is not really a serious one. Our method is here less elaborate than the one of
\cite{CCH}, its give a more general result but it is less sharp :
there are two restrictions : the first one is the dimension restriction $\nu>3$ which is unsatisfactory and the
second concerns the limitation  $p<\nu$ which is perhaps also  unsatisfactory when 
$M$ has only one end. However there are recent results of T. Coulhon and N. Dungey in this direction
\cite{CDu}. 

There is now a long list of complete Riemannian manifolds $(M_0,g_0)$ on which the Riesz transform
is bounded on $L^p$ for every $p\in ]1,\infty[$ and satisfying our hypothesis. For instance  Cartan-Hadamard manifolds with a spectral gap \cite{Lo},
 non-compact symmetric spaces \cite{Anker}  and Lie 
 groups of polynomial growth \cite{Alex},
manifolds with nonnegative Ricci curvature and maximal volume growth \cite{Bakry} (see the discussion at the end of the
proof of theorem \ref{main} about the case of manifolds with nonnegative Ricci curvature 
and non maximal volume growth).
Also H.-Q. Li \cite{Li} proved that the Riesz transform on $n$-dimensional cones with compact basis is bounded on $L^p$ for $p < p_0$, where $$p_0 =  \begin{cases} n \Big( \frac{n}{2} - \sqrt{\big(\frac{n-2}{2} 
    \big)^2 + \lambda_1} \Big)^{-1}, \quad \lambda_1 < n-1 \\+\infty, \qquad \lambda_1 \geq n-1, \end{cases}$$
    where $\lambda_1$
     is the smallest nonzero eigenvalue of the Laplacian on the basis. Note that  $p_0>n$.
Our proof also applies to manifold isometric at infinity to several copies of cones, hence our theorem \ref{main} also gives a partial answer to the open problem 8.1 of \cite{CCH} :
\begin{cor} 
  If $(M,g)$ is a  smooth Riemannian manifold of dimension $n\ge 4$ with conic ends, then the Riesz transform is bounded
  on $L^p$ for any $p\in ]1,n[$.
\end{cor}

\section{Analytic preliminaries}
\subsection{A Sobolev inequality}
\begin{prop}\label{sob} Let $(M,g)$ be a complete Riemannian manifold with Ricci curvature bounded from below
and with positive injectivity radius, then for any $p\in [1,\infty]$, there is a constant $C$ such that
for all $\varphi \in C^\infty_0(M)$
$$\|df\|_{L^p}\le C\left[\|\De f\|_{L^p}+\|f\|_{L^p}\right].$$
\end{prop}
\begin{rem}
\begin{enumeroman}
\item T. Coulhon and X. Duong have shown that for every complete Riemannian manifolds and any $p\in ]1,2]$,
there is a constant $C$ such that 
$$\forall f\in \Ci, \|df\|_{L^p}^2\le C\|\De f\|_{L^p}\|f\|_{L^p}.$$
When $p\in]1,2]$, this is clearly a stronger result. 
\item We don't know whether our geometric condition are necessary.
\item This lemma belongs to the folklore, see for instance \cite{H} for the same idea.
 E.B. Davies has proven this result for manifold with bounded geometry (see corollary 10 in \cite{Da}).
\end{enumeroman}
\end{rem}

\proof According to the results of Anderson-Cheeger \cite{AC}, there is a $r_H>0$  depending only on the lower bound of
the Ricci curvature and the injectivity radius such that for any $x\in M$ there is a harmonic coordinate
chart on  the geodesic ball 
$B(x,r_H)$ :
$$H \,:\, B(x,r_H)\rightarrow \R^n$$ with $H(x)=0,$
 $$\|H_*g-\eucl\|_{C^{\frac12}}\le C$$ and
 $$\frac14 \eucl\le H_*g\le 4\, \eucl,$$
 where $\eucl$ denotes the Euclidean metric on $\R^n$.
Let $\Delta_0=-\sum_{i}\frac{\partial^2}{\partial x_i^2}$ the Euclidean Laplacian and
$$\tilde\Delta=-\sum_{i,j} \tilde{g}^{i,j}\frac{\partial^2}{\partial x_i\partial x_j}$$
be the Laplacian of the metric $\tilde g=H_*g$ on $H(B(x,r_H))$.
We define $U_r=H(B(x,r))$ and $V_r=B(0,r)$ the Euclidean ball of radius $r$, we have
$$V_{r/2}\subset U_r\subset V_r.$$
Moreover, on $V_r$ we always have the Sobolev inequality
$$\forall f\in C_0^\infty(V_r), 
 r^2\|\partial^2 f\|_{L^p}+r\|\partial f\|_{L^p}+\| f\|_{L^p}\le C(n,p)\left[ r^2\|\De_0f\|_{L^p}+\|f\|_{L^p}\right].$$
But by hypothesis  if $r<r_H$ and $f \in C_0^\infty(V_r)$, we have
$$\|\De_0f\|_{L^p}\le\|\tilde\De f\|_{L^p} + \|(\De_0-\tilde\De )f\|_{L^p}\le\|\tilde\De f\|_{L^p} +C\sqrt{r}
 \|\partial^2 f\|_{L^p} $$
 Hence there is a $r_0>0$ depending only on $p,n$, on the lower bound of
the Ricci curvature and on the injectivity radius such that
for any $r<r_0$ 
$$\forall f\in C_0^\infty(V_r):\ r^2\|\partial^2 f\|_{L^p}+r\|\partial f\|_{L^p}+\| f\|_{L^p}\le C\left[ r^2
\|\tilde \De f\|_{L^p}+\|f\|_{L^p}\right].$$
But on $V_{r_0}$ the metric $\tilde g$ and $\eucl$ are quasi isometric, hence
 $$\forall f\in C_0^\infty(V_r):r \|df\|_{L^p}\le C\left[ r^2
\|\tilde \De f\|_{L^p}+\|f\|_{L^p}\right].$$
where all the norms are now measured with respect to the norms associated to the metric $\tilde g$.
Coming back to $(M,g)$, we have proved that there is constants $\rho>0$ and $C$ such that
for any $x\in M$, 
$$\forall f\in C_0^\infty(B(x,\rho)),\ \|d f\|_{L^p(B(x,\rho))}\le C\left[ 
\|\De f\|_{L^p(B(x,4\rho))}+\|f\|_{L^p(B(x,4\rho))}\right].$$

Now a classical covering argument using the Bishop-Gromov estimate implies the desired result.
 \endproof

\subsection{Some estimate on the Poisson operator}
\begin{lem}\label{PE1}
Let $(M,g)$ be a complete Riemannian manifold which for some $\nu >2$ and $C>0$ satisfies the Sobolev
inequality :
$$\forall \varphi\in  \Ci, C\|\varphi\|_{L^{\frac{2\nu}{\nu-2}}}\le \|d\varphi\|_{L^2}$$
then the Schwarz kernel $P_\sigma(x,y)$ of the Poisson operator $e^{-\sigma\sqrt{\De }}$ satisfies
$$P_\sigma(x,y)\le \frac {C\sigma}{\left( \sigma^2+d(x,y)^2 \right)^{\frac{\nu+1}{2}}}.$$
Moreover if $1\le r\le p\le +\infty$ then 
$$\left\|e^{-\sigma\sqrt{\De }} \right\|_{L^r\to L^p}\le \frac{C}{\sigma^{\nu\left(\frac1r-\frac1p\right)}}.$$
\end{lem}
We know that the heat operator $e^{-t\De}$ and the Poisson operator are related through the
subordination identity :
$$ e^{-\sigma\sqrt{\De }}=\frac{\sigma}{2\sqrt\pi}\int_0^\infty e^{-\frac{\sigma^2}{4t}} e^{-t\De}
\frac{dt}{t^{3/2}}.$$
Hence these properties follow directly from the corresponding ones for the heat operator $e^{-t\De}$ and its
Schwarz kernel $H_t(x,y)$ :
$$H_t(x,y)\le \frac{c}{t^{\nu/2}}e^{-\frac{d(x,y)^2}{5t}}$$ and
if  $1\le r\le p\le +\infty$ then 
$$\left\|e^{-t{\De }} \right\|_{L^r\to L^p}\le \frac{C}{t^{\frac{\nu}{2}\left(\frac1r-\frac1p\right)}},$$
which are consequences of the Sobolev inequality \cite{V, CSV}. 

We will also need an estimate for the derivative of the Poisson kernel :
\begin{lem}\label{PE2}
Under the assumptions of lemma (\ref{PE2}), let $\Omega\subset M$ be a open subset and
$K$ be a compact set in in the interior of $M\setminus\Omega$ then
$$\left\|e^{-\sigma\sqrt{\De }} \right\|_{L^p(\Omega)\to L^\infty(K)}\le \frac{C}{(1+\sigma)^{\nu/p}},$$
$$\left\|\nabla e^{-\sigma\sqrt{\De }} \right\|_{L^p(\Omega)\to L^\infty(K)}\le \frac{C}{(1+\sigma)^{\nu/p}}.$$
\end{lem}
\proof 
The first identity is only a consequence of the first lemma. To prove the second inequality, we will again
only show the corresponding estimate for the heat operator.
First, according the local Harnack inequality (see V.4.2 in \cite{CSV}), there is a constant $C$ such that
for any $x\in K$, $t\in]0,1]$ and $y\in M$ :
$$\left|\nabla_x H_t(x,y)\right|\le \frac{C}{\sqrt{t}} H_{2t}(x,y).$$
But by assumption there is a constant $\eps>0$ such that
for all $(x,y)\in K\times\Omega$ then
$$d(x,y)\ge \eps$$ hence
for all $(x,y)\in K\times\Omega$ then
$$H_{2t}(x,y)\le  \frac{c}{t^{\nu/2}}e^{-\frac{\eps^2}{10t}},$$
we easily obtain that there is a certain constant $C$ such that 
$$\forall t\in ]0,1] : \left\|\nabla e^{-t{\De }} \right\|_{L^p(\Omega)\to L^\infty(K)}\le C$$
Now assume that $t>1$ :
$$\left\|\nabla e^{-t{\De }} \right\|_{L^p(\Omega)\to L^\infty(K)}\le
\left\|\nabla e^{-\frac12{\De }} \right\|_{L^\infty(M)\to L^\infty(K)} 
\left\| e^{-(t-\frac12){\De }} \right\|_{L^p(\Omega)\to L^\infty(M)}.$$
But we have
$$\left\| e^{-(t-\frac12){\De }} \right\|_{L^p(\Omega)\to L^\infty(M)}\le \frac{c}{(t-1/2)^{\nu/p}}.$$
And
because
$$\sup_{x\in K} \int_M |\nabla_x H_{1/2}(x,y)| dy\le C \sup_{x\in K} \int_M H_{1}(x,y) dy\le C,$$ we get the
desired result.
\endproof
\section{Proof of the main theorem }
Let $(M_0,g_0)$ be a complete Riemannian manifold, we assume that the Ricci curvature of $(M_0,g_0)$ is bounded from
below , that the injectivity radius of $(M_0,g_0)$ is positive and that for some $\nu>3$ and $C>0$, that $(M,g)$ satisfies the Sobolev inequality
$$\forall \varphi\in  C^\infty_0(M_0),\  C\|\varphi\|_{L^{\frac{2\nu}{\nu-2}}}\le \|d\varphi\|_{L^2}.$$
We assume that on $(M_0,g_0)$ the Riesz transform is bounded on $L^p$ for some $p\in ]\nu/(\nu-1), \nu[$.
And we consider $M$ a complete Riemannian manifold such that
outside  compact sets $K\subset M$ and $K_0\subset M_0$, 
$M\setminus K$ is isometric to $k$ copies of $M_0\setminus K_0$. 
We are going to prove that on $M$ the Riesz
transform is also bounded on $L^p$. The first step is to build a good parametrix for the Poisson operator on
$M$.
The first problem is that the operator $\sqrt{\Delta}$ is not a differential operator, we circonvent this
difficulties by working on $\R_+\times M$.  As a matter of fact, the Poisson operator solves
the Dirichlet problem :
$$\left\{\begin{array}{ll}
\left(-\frac{\partial^2}{\partial \sigma^2}+\De\right) u(\sigma,x)=0& \mbox{on } ]0,\infty [ \times M\\
u(0,x)=u(x)& \\
\lim_{\sigma\to\infty} u(\sigma,.)=0& \\
\end{array}\right.$$
The construction of the parametric will be standard, the nature of the operator $\sqrt{\Delta}$
implies that we can not used the Duhamel formula, instead we used the Green operator.
The idea is to find an approximate solution for this problem $E(u)$ and
then to used the fact that if
$G$ is the Green operator of the operator $-\frac{\partial^2}{\partial \sigma^2}+\De$
for the Dirichlet boundary condition then
$$e^{-\sigma\sqrt{\De}}u=E(u)+G( -\frac{\partial^2}{\partial \sigma^2}+\De)E(u).$$

\subsection{The parametrix construction}
Let $\tilde K$ be a another compact set in $M$ containing $K$ in its interior.
Let $E_1,...,E_b\subset M$ be $k$ pieces of $M\setminus K$, each one being isometric to $M_0\setminus K_0$.
Let $\rho_0,...,\rho_1$ a smooth partition of unity such that
$$\supp \rho_0\subset \tilde K\mbox{ and }
\forall i\ge 1,\ \supp\rho_i\subset E_i$$
Let also $\varphi_0,...,\varphi_1$ be smooth function, such that
$$\supp \varphi_0\subset \tilde K \mbox{ and }
\forall i\ge 1,\ \supp\varphi_i\subset E_i$$
We moreover require that $$\varphi_i\rho_i=\rho_i.$$
Let $\Delta_0$ be the realization of the Laplace operator on $\tilde K$ for the Dirichlet boundary
condition. And let $\Delta_i$ be the Laplace operator on the $i^{th}$ copy of $M_0$.  Let 
$e^{-\sigma \sqrt{\De_i}}$ its associated Poisson operator then we 
define for $u\in L^p(M)$:
$$E(u)(\sigma)=\sum_{i=0}^b \varphi_i (e^{-\sigma \sqrt{\De_i}}\rho_i u).$$
We can easily computed :
$$\left(-\frac{\partial^2}{\partial \sigma^2}+\De\right)E(u)=\sum_{i=0}^b [\Delta,\varphi_i](e^{-\sigma
\sqrt{\De_i}}\rho_i u)=f(\sigma,x)=\sum_{i=0}^b f_i(\sigma,x),$$
where
$$f_i(\sigma,x)=[\Delta,\varphi_i](e^{-\sigma
\sqrt{\De_i}}\rho_i u)(x)=\Delta\varphi_i (x)(e^{-\sigma
\sqrt{\De_i}}\rho_i u)(x)-2\left\langle d\varphi_i(x),d(\sqrt{\De_i}\rho_i u)(x)\right\rangle.$$
From the lemma \ref{PE2}, we easily get for $i\ge 1$ and 
all $\sigma\ge 0 $:
\begin{equation}\label{esti1}
\|f_i(\sigma)\|_{L^1}+\|f_i(\sigma)\|_{L^p}\le \frac{C}{(1+\sigma)^{\nu/p}} \|\rho_i u\|_{L^p}.
\end{equation}

\noindent Let's us explain why this estimate also hold for $f_0$. Note that the operator 
$$\cS(\sigma)=[\Delta,\varphi]e^{\sigma\sqrt{\De_0}}\rho_0$$ is
a operator with smooth Schwarz kernel and proper support, moreover 
it is bounded when $\sigma \to 0$.
Hence there is a constant $C$ such that 
$$\forall \sigma\in [0,1], \|\cS(\sigma)\rho_0 u\|_{L^\infty}\le C\|\rho_0 u\|_{L^p}.$$
Now the operator $\De_0$ has a spectral gap on $L^p$ 
(its $L^p$ spectrum is also its $L^2$ spectrum), hence
there is a constant
$C$ such that for all
$\sigma\ge 0$ then
$$\|e^{-\sigma\sqrt{\De_0}}\|_{L^p\to L^p} \le C e^{-\sigma/C},$$

\noindent Hence for $\sigma \ge 1$ :
$$\| \cS(\sigma)u\|_{L^\infty}\le \|[\Delta,\varphi]e^{\frac12 \sqrt{\De_0}}\|_{L^p\to L^\infty} 
\|e^{-(\sigma-1/2)\sqrt{\De_0}}\rho_0u\|_{L^p}\le C e^{-\sigma/C} \|\rho_0 u\|_{L^p}.$$
The result follows by noticing that 
 the $f_i$ have compact support in $\tilde K\setminus K$. Eventually we obtain the estimate :
\begin{lem}
When $u\in L^p(M)$ and let $S(u)=f$ then 
$$\forall \sigma\ge 0,\ \|S(u)\|_{L^1}+\|S(u)\|_{L^p}\le \frac{C}{(1+\sigma)^{\nu/p}} \|u\|_{L^p}.$$
\end{lem}
\subsection{The Riesz transform on $M$}
We introduce now $G$ the Green operator on $\R_+\times M$ for the Dirichlet boundary condition.
Its Schwarz kernel is given by 
$$G(\sigma,s,x,y)=\int_ 0^\infty \left[\frac{e^{-\frac{(\sigma-s)^2}{4t}} -e^{-\frac{(\sigma+s)^2}{4t}}}{\sqrt{4\pi t}}\right] H_t(x,y) dt$$
where 
$H_t$ is the heat kernel on $M$ and 
$$\frac{e^{-\frac{(\sigma-s)^2}{4t}}-e^{-\frac{(\sigma+s)^2}{4t}}}{\sqrt{4\pi t}}$$ the heat kernel on the half-line
$\R_+$ for the Dirichlet boundary condition.
 We have
 $$ e^{-\sigma\sqrt{\De}}u=E(u)(\sigma)+G(S(u)).$$
 Hence
 $$\De^{-1/2}u=\int_0^\infty e^{-\sigma\sqrt{\De}}ud\sigma=
 \sum_{i=0}^b \varphi_i \De_i^{-1/2}\rho_i u+\int_{\R^2_+\times M} G(\sigma,s,x,y)f(s,y)d\sigma dsdy.$$
But \begin{equation*}
\begin{split}
\int_0^\infty G(\sigma,s,x,y)d\sigma&=\frac{1}{\sqrt{4\pi}}\int_0^\infty  \left[\int_{-s}^s
e^{-\frac{v^2}{4t}}dv\right] H_t(x,y) \frac{dt}{\sqrt{t}}\\
&=\frac{2}{\sqrt{\pi}}\int_0^\infty e^{-r^2}\left[ \int_{0}^{\frac{s^2}{4r^2}}H_t(x,y) dt\right] dr\\
\end{split}\end{equation*}
Let 
\begin{equation*}
\begin{split}g(x)&=\int_{\R^2_+\times M} G(\sigma,s,x,y)f(s,y)d\sigma dsdy\\
&=\frac{2}{\sqrt{\pi}}\int_{\R^2_+} e^{-r^2}\left[ \int_{0}^{\frac{s^2}{4r^2}}(e^{-t\De}f(s))(x)dt\right]
drds,\\
\end{split}\end{equation*}
so that 
$$\De^{-1/2} u=\sum_{i=0}^b \varphi_i \De_i^{-1/2}\rho_i u+g.$$
The following lemma is now the last crucial estimate:
\begin{lem}\label{lastest}There is a constant $C$ such that 
$$\|\De g\|_{L^p}+\| g\|_{L^p}\le C\|u\|_{L^p}.$$
\end{lem}
\proof Recall that according to \cite{carronduke}, $(M,g)$ itself satisfies the same Sobolev inequality~:
$$\forall \varphi\in  \Ci,\  C\|\varphi\|_{L^{\frac{2\nu}{\nu-2}}}\le \|d\varphi\|_{L^2}.$$
Hence the heat operator satisfies the following mapping properties :
for $1\le q\le p\le +\infty$ we have
$$\left\|e^{-t{\De }} \right\|_{L^q\to L^p}\le \frac{C}{t^{\frac{\nu}{2}\left(\frac1q-\frac1p\right)}}.$$
As a consequence,
for all $t\in [0,1]$, then 
$$\|(e^{-t\De}f(s))\|_{L^p}\le \|f(s))\|_{L^p}\le \frac{C}{(1+s)^{\nu/p}}\,\|u\|_{L^p}$$
and if $t>1$, then
$$\|(e^{-t\De}f(s))\|_{L^p}\le \left\|e^{-t\De}\right\|_{L^1\to L^p}\|f(s))\|_{L^1}\le
 \frac{1}{t^{\frac{\nu}{2}(1-\frac1p)}}\frac{C}{(1+s)^{\nu/p}}\,\|u\|_{L^p},$$
 Hence 
\begin{equation*}
\begin{split}
\|g\|_{L^p}&
\le \frac{2}{\sqrt{\pi}}\int_{\R^2_+} e^{-r^2}\left[
\int_{0}^{\frac{s^2}{4r^2}}\|(e^{-t\De}f(s))\|_{L^p}dt\right] dsdr\\
&\le  \frac{2}{\sqrt{\pi}}\int_0^ \infty e^{-r^2}\left[
\int_{0}^{\frac{s^2}{4r^2}}
\frac{C}{\max \left(1, t^{\frac{\nu}{2}(1-\frac1p)}\right)}\frac{1}{(1+s)^{\nu/p}}dt\right]
 dsdr\, \|u\|_{L^p} .\\
\end{split}\end{equation*}
But because $p<\nu$, we have
\begin{equation*}
\begin{split}
&\int_{ \{2r\sqrt{t}\le s\} } e^{-r^2}\frac{1}{\max \left(1,
t^{\frac{\nu}{2}(1-\frac1p)}\right)}\frac{1}{(1+s)^{\nu/p}}dsdtdr \\
&=\frac{\nu}{\nu-p} \int_{\R_+^2} e^{-r^2} \frac{1}{\max \left(1,
t^{\frac{\nu}{2}(1-\frac1p)}\right)}\frac{1}{(1+2r\sqrt{t})^{\nu/p-1}}dtdr\\
 \end{split}\end{equation*}
and this integral is finite exactly when $p>\nu/(\nu-1)$ and $\nu>3$. 

It remains to estimate $\|\De g\|_{L^p}$, this is easier because 
\begin{equation*}
\begin{split}\De g&=\frac{2}{\sqrt{\pi}}
\int_{\R^2_+} e^{-r^2}\left[ \int_{0}^{\frac{s^2}{4r^2}}\De (e^{-t\De}f(s))dt\right] drds\\
&=-\frac{2}{\sqrt{\pi}}
\int_{\R^2_+} e^{-r^2}\left[ \int_{0}^{\frac{s^2}{4r^2}}\frac{d}{dt}(e^{-t\De}f(s))dt\right] drds\\
&=\frac{2}{\sqrt{\pi}}
\int_{\R^2_+} e^{-r^2}\left[f(s)- (e^{-\frac{s^2}{4r^2}\De}f(s))\right] drds
\end{split}\end{equation*}
Hence 
\begin{equation*}
\begin{split}\|\De g\|_{L^p}&\le\frac{4}{\sqrt{\pi}}
\int_{\R^2_+} e^{-r^2}\|f(s)\|_{L^p}  drds\\
&\le  \frac{4}{\sqrt{\pi}}
\int_{\R^2_+} e^{-r^2}\frac{C}{(1+s)^{\nu/p}}drds\ \|u\|_{L^p} .\\
\end{split}\end{equation*}
 \endproof
Now we can finish the proof of the main theorem :
Let $T_i$ be the Riesz transform associated with the operator $\De_i$ then we obtain
$$d\De^{-1/2}u=\sum_{i=0}^b \varphi_i T_i\rho_i u+\sum_{i=0}^bd\varphi_i(\De^{-1/2}_i\rho_i u)+dg.$$
By hypothesis, $T_i$ is bounded on $L^p$ for all $i\ge 1$. Moreover $\varphi_0 T_i\rho_0$
is a pseudo differential operator of order $0$ with proper support it is also bounded on $L^p$.
Moreover, the Sobolev inequality 
$$\forall \varphi\in  \Ci,\  C\|\varphi\|_{L^{\frac{2\nu}{\nu-2}}}\le \|d\varphi\|_{L^2}.$$
also implies the following mapping properties of the $\De_i^{-1/2}$ (\cite{V}) :
$$\left\| \De^{-1/2}_i\right\|_{L^p\to L^{\frac{p\nu}{\nu-p}}}\le C.$$
Hence 
$$\left\| d\varphi_i(\De^{-1/2}_i\rho_i u)\right\|_{L^p}\le C \|\De^{-1/2}_i\rho_i u\|_{L^p(\tilde K)}
\le C' \|\De^{-1/2}_i\rho_i u)\|_{L^{\frac{p\nu}{\nu-p}}(\tilde K)}\le C \|\rho_i u\|_{L^p}.$$
Moreover the lemmas (\ref{lastest}) and (\ref{sob}) implie that
$$\|dg\|_{L^p}\le C \|u\|_{L^p}$$
All these estimates yield the fact that the Riesz transform is bounded on $L^p$.
\subsection{A comment on manifold with non negative Ricci curvature}
The proof of theorem C is fairly general, we can easily make a list of the properties which make it runs ;
let $(M_i,g_i)$ $i=1,...,b$ be complete Riemannian manifolds and let $(M,g)$ be isometric at infinity to the
disjoint union $M_1\cup...\cup M_b$. That is to say there are compact sets $K\subset M$, $K_i\subset M_i$
such that $M\setminus K$ is isometric to $(M_1\setminus K_1)\cup...\cup (M_b\setminus K_b)$.
Let $\tilde K\subset \widehat K$ such that $\tilde K$ (resp. $\widehat K$) contains $K$ in its interior
(resp. $\tilde K$). And let $\widehat K_i,\tilde K_i\subset M_i$ such that :
$$M\setminus \tilde K\simeq (M_1\setminus\tilde K_1)\cup...\cup (M_b\setminus\tilde K_b),\ M\setminus\widehat
K\simeq (M_1\setminus\widehat K_1)\cup...\cup (M_b\setminus\widehat K_b),$$ let $\Delta_i$ be the Laplace
operator on $M_i$.
 We assume that on each $M_i$, the Ricci curvature is bounded from below and that the injectivity radius is
 positive, such that on each $M_i$ and $M$, we get the estimate induced by the 
Sobolev inequality (\ref{sob}).
Assume that for some functions $f,g\,:\,\R_+\rightarrow \R_+^*$ we have the estimate :

$$\left\|e^{-\sigma \sqrt{\De_i}}\right\|_{L^p(M_i\setminus \widehat K_i)\to L^\infty(\tilde K_i)}
+\left\|\nabla e^{-\sigma \sqrt{\De_i}}\right\|_{L^p(M_i\setminus \widehat K_i)\to L^\infty(\tilde K_i)}
\le \frac{1}{f(\sigma)}\ ,$$
and that on the manifold $M$ :
$$\left\|e^{-t \De}\right\|_{L^1(\widehat K)\to L^p(M)}\le \frac{1}{g(t)}.$$
with
\begin{equation}\label{intcondi1}
\int_0^\infty \frac{ds}{f(s)}<\infty\end{equation}
\begin{equation}\label{intcondi2} 
\int_{\R^2_+} e^{-u^2} \min\left(1,\frac{1}{g(t)}\right)\left[\int_{2u\sqrt{t}}^\infty \frac{ds}{f(s)}\right]
dudt<\infty.\end{equation}

Then if for all $i$, the Riesz transform $T_i:=d\De_i^{-1/2}$ is bounded on $L^p$, then on $M$, the Riesz
transform is also bounded on $L^p$.

A natural and well study class of manifolds satisfying such estimates are manifolds satisfying the
so called relative Faber-Krahn inequality: for some $\alpha>0$ and $C>0$, we have :
 \begin{equation*}
\forall B(x,R),\forall \Omega\subset B(x,R),\ 
\lambda_1(\Omega)\ge \frac{c}{R^2} \left(\frac{\vol \Omega}{\vol B(x,R)}\right)^{-\alpha}
\end{equation*}
where $$\lambda_1(\Omega)=\inf_{f\in C^\infty_0(\Omega)} \frac{\int_\Omega|df|^2}{\int_\Omega f^2} $$
is the first eigenvalue of the Laplace operator on $\Omega$ for the Dirichlet boundary condition.
According to A. Grigor'yan \cite{G} this inequality is equivalent to the conjunction of the doubling
property: uniformly in $x$ and $R>0$ we have
$$\ \frac{\vol B(x,2R)}{\vol B(x,R)}\le C$$
and of the upper bound on the heat operator
$$H_t(x,y)\le \frac{C}{\vol B(x,\sqrt{t})}e^{-\frac{d(x,y)^2}{5t}}\ .$$
Manifolds with non negative Ricci curvature are examples on manifolds satisfying this 
relative Faber-Krahn inequalities. 

Assume that each $M_i$ satisfies this relative Faber-Krahn inequalities and if we assume
that for $i=1,...,b$, there is a point $o_i\in K_i$ and all $R\ge 1$ 
$$\vol B(o_i,R):=V_i(R)\ge C R^\nu$$ 
then we get easily from the subordination identity :
$$\left\|e^{-\sigma \sqrt{\De_i}}\right\|_{L^p(M_i\setminus \widehat K_i)\to L^\infty(\tilde K_i)}
+\left\|\nabla e^{-\sigma \sqrt{\De_i}}\right\|_{L^p(M_i\setminus \widehat K_i)\to L^\infty(\tilde K_i)}
\le \frac{1}{(1+\sigma)^{\nu/p}}.$$
Now the problem comes from the fact that we don't know how to obtain a relative Faber-Krahn inequality on $M$
from the one we assume on the $M_i$'s. However, recently in \cite{GS2}, A. Grigor'yan and L. Saloff-Coste 
have announced the following very useful result (see also\cite{GS}) : when the $M_i$'s satisfy the 
relative Faber-Krahn inequality then 
\begin{equation*}
\forall B(x,R)\subset M,\forall \Omega\subset B(x,R),\ \ 
\lambda_1(\Omega)\ge \frac{c}{R^2} \left(\frac{\vol \Omega}{\mu(x,R)}\right)^{-\alpha}
\end{equation*}
Where $$\mu(x,R)=\left\{ \begin{array}{ll}
\vol B(x,R)&\mbox{if }B(x,R)\subset M\setminus K\\
\inf_i V_i(R)&\mbox{else }  \\
\end{array}\right. $$ Hence 
from our volume growth estimate, we will obtain (see \cite{G}) when $t\ge 1$ :
$$\left\|e^{-t \De}\right\|_{L^1(\widehat K)\to L^p(M)}\le
\frac{C}{t^{\frac{\nu}{2}\frac{p-1}{p}}}$$
With this result of A. Grigor'yan and L. Saloff-Coste and with the result of D. Bakry
\cite{Bakry}, we will obtain :
\begin{prop} Let $(M_1,g_1),..., (M,g_b)$ be complete Riemannian manifolds with non negative Ricci curvature.
Assume that on all $M_i$ we have the volume growth lower bound :
$$\vol B(o_i,R)\ge C R^\nu.$$
Then assume that  $\nu>3$ then on any manifold isometric at infinity to the disjoint union
of the $M_i$'s, the Riesz transform is bounded on $L^p$ for all $p\in ]\nu/(\nu-1),\nu[$.
\end{prop}
Note that we have remove the hypothesis on the injectivity radius : as matter of fact the Sobolev inequality
(\ref{sob}) holds on manifolds with non negative Ricci curvature (\cite{CDfull}) ; then it is easy to show
that this inequality also hold on $M$.

\end{document}